\documentclass[letterpaper, 10 pt, conference]{ieeeconf}  
\IEEEoverridecommandlockouts        
\usepackage{graphicx} 
\usepackage{epstopdf} 
\usepackage{psfrag} 
\usepackage{amsmath} 
\usepackage{amssymb}  
\usepackage{amsthm,thmtools,dsfont}
\usepackage{color}
\usepackage{tabu}
\usepackage[hidelinks]{hyperref}
\usepackage{cleveref}
\let\labelindent\relax
\usepackage{enumitem} \setlist[itemize,1]{leftmargin=\dimexpr 24pt}
\usepackage{cite}
\usepackage{microtype,enumerate,url}
\usepackage{BOONDOX-ds}
\usepackage{tikz} 
\usetikzlibrary{patterns} 
\usetikzlibrary{positioning} 
\usepackage[fulladjust]{marginnote}
\usepackage{pgfplots} %
\pgfplotsset{compat=1.17} 
\usepackage{pgfplotstable} 
\usepackage{xstring} 
\usepackage{booktabs} 
\usepackage{colortbl}
\usepackage[skip=0pt]{caption}
\usepackage{mathdots} 
\usepackage{tcolorbox}
\usepackage[normalem]{ulem} 
\allowdisplaybreaks 
\usepackage{tikz-cd}
\usepackage[switch, modulo]{lineno}

\graphicspath{{./imgs/}}
\usepackage{algorithm}
\usepackage{algorithmic}

\renewcommand{\theenumi}{(\roman{enumi})}

\makeatletter
\def\@IEEEsectpunct{.\ \,}
\def\paragraph{\@startsection{paragraph}{4}{\z@}{1.5ex plus 1.5ex minus 0.5ex}%
{0ex}{\normalfont\normalsize\itshape}}
\makeatother

\renewcommand*\descriptionlabel[1]{\hspace\labelsep\normalfont#1}

\declaretheorem[style=definition]{theorem}
\declaretheorem[style=definition]{corollary}

\declaretheorem[style=definition]{lemma}
\declaretheorem[style=definition,qed=$\vartriangle$]{remark}

\declaretheorem[style=definition,numbered=no]{standing assumption}

\declaretheorem[style=definition]{problem}

\renewcommand\thmcontinues[1]{continued}

\newcommand {\nn}{\nonumber}
\newcommand{\beq}{\begin{equation}}
\newcommand{\eeq}{\end{equation}}
\newcommand {\bseq}{\begin{subequations}}
\newcommand {\eseq}{\end{subequations}}
\newcommand {\bma}{\left[}
\newcommand {\ema}{\right]}

\newcommand {\N}{\mathbb{N}} 	
\newcommand {\R}{\mathbb{R}} 	

\newcommand {\T}{\mathbb{T}} 	
\newcommand {\W}{\mathbb{W}} 	
\newcommand {\B}{\mathcal{B}} 	
\newcommand {\BL}{\B|_L}
\renewcommand{\L}{\mathfrak{L}} 	

\newcommand{\Reference}{\mathcal{R}} 
\newcommand{\Plant}{\mathcal{P}}
\newcommand{\Controller}{\mathcal{C}} 
\newcommand{\Hidden}{\mathcal{N}} 
\newcommand{\PC}{\mathcal{P} \lVert_{c} \mathcal{C}} 
\newcommand{\PR}{\mathcal{P} \lVert_{w} \mathcal{R}} 
\newcommand{\PlantL}{{\Plant|_L}}
\newcommand{\HiddenL}{{\Hidden|_L}}
\newcommand{\ReferenceL}{{\Reference|_L}}
\newcommand{\ControllerL}{{\Controller|_L}}

\renewcommand{\implies}{\Rightarrow} 	
\newcommand{\Bell}{\text{\boldmath$\ell$\unboldmath}} 
\newcommand{\Image}{\operatorname{im}} 
\newcommand{\rank}{\operatorname{rank}} 
\newcommand{\blockdiag}{\operatorname{block-diag}} 

\title{\Large {\bfseries  Controller implementability: a data-driven approach}}

\author{Alberto Padoan,  Jeremy Coulson, and Florian D\"orfler\thanks{
		A.  Padoan and  F.  D\"orfler are with 
		the Department of Information Technology and Electrical Engineering at		
		ETH Z\"urich,  Z\"urich, Switzerland 
		{\tt\footnotesize  \{apadoan, dorfler\}@control.ee.ethz.ch}. 
		J.  Coulson, 
		is with the Department of Electrical and Computer Engineering at the University of Wisconsin-Madison, USA.
		{\tt\footnotesize  jeremy.coulson@wisc.edu}. 
  }
}
\date{\small\today} 

\begin{document}

\maketitle
\thispagestyle{empty}
\pagestyle{empty}

\begin{abstract}
We study the controller implementability problem, which seeks to determine if a controller can make the closed-loop behavior of a given plant match that of a desired reference behavior.
 We establish necessary and sufficient conditions for controller implementability which only rely on raw data.  Subsequently, we consider the problem of constructing controllers directly from data. By leveraging the concept of canonical controller, we provide a formula to directly construct controllers that implement plant-compatible reference behaviors using measurements of both reference and plant behaviors.
\end{abstract}

\section{Introduction}
The problem of control design can be split into three parts: (i) to describe the set of admissible controllers; (ii) to describe the properties that the controlled system should have; and (iii) to find an admissible controller such that the resulting controlled behavior has desired properties~\cite{willems1997interconnections}. 
Generally, in order to solve a control design problem, one needs access to a model of the system to be controlled and a model of the reference behavior. 
However, in many  situations of practical interest  obtaining such models is expensive, time-consuming, or simply impossible, and the control designer only has access to measured data~\cite{markovsky2021behavioral}. 
This has motivated the development of new direct data-driven control methods that bypass system identification and aim to compute controllers directly from data, see, \textit{e.g.}, the recent survey~\cite{markovsky2021behavioral}.

The paper studies the controller implementability problem~\cite{willems2002synthesis}. The problem is to find, if possible, a controller which makes the closed-loop behavior of a plant equal to that of a desired reference behavior. While the problem does admit a model-based solution, our objective is to provide an alternative solution that is compatible with modern data-driven approaches.
Using the language of behavioral systems theory~\cite{willems1986timeI}, we regard 
finite-horizon behaviors of finite-dimensional, linear, time-invariant (LTI) systems as subspaces represented by raw data matrices~\cite{markovsky2021behavioral}. We establish necessary and sufficient conditions for controller implementability which can be tested directly from raw data. Furthermore, we also consider the problem of constructing controllers directly from data.  We provide a formula to directly construct controllers that implement plant-compatible reference behaviors using only measurements of the reference and plant behaviors.

\textbf{Contributions}: 
The contributions of the paper are twofold. We establish new necessary and sufficient conditions for solving the controller implementability problem, thus characterizing all implementable controlled behaviors in both model-based and data-driven scenarios. We provide a formula for a canonical controller that implements any given reference behavior, whenever this is possible; the controller depends solely on the reference and plant behaviors and can be directly obtained from data. 

\textbf{Related work}: 
The controller implementability problem has been originally studied in~\cite{willems2002synthesis}, where necessary and sufficient conditions for implementability are given for continuous-time behaviors. 
The concept of a canonical controller has been implicitly defined in the seminal paper~\cite{willems2002synthesis} and, subsequently, formalized for general systems, e.g., in~\cite{van2003achievable} and~\cite{julius2005canonical}. Our results extend existing results presented in~\cite{willems2002synthesis, van2003achievable, julius2005canonical}
to discrete-time, finite-horizon behaviors, eliminating the need for parametric models and enabling the direct use of raw data.  Over the past two decades, the data-driven approaches have received increasing attention, primarily due to the surge in availability of data, see, \textit{e.g.} the recent survey~\cite{markovsky2021behavioral}.
A simple, yet paradigmatic instance of the controller implementability problem is the \textit{exact model matching} problem~\cite{wolovich1972use,wang1972exact}, whereby one seeks a state feedback law for a given finite-dimensional LTI system to make the closed-loop transfer function equal to a given transfer function. The problem has been widely studied in a model-based context~\cite{wolovich1972use,wang1972exact,wang1972solution} and  it is well-known that the problem can be reduced to solving a set of linear algebraic equations~\cite{wang1972exact,wang1972solution}. 
The recent paper~\cite{breschi2021direct} presents analogous findings in a data-driven context. Our results generalize the findings of~\cite{breschi2021direct} in a representation-free setting, without requiring plant and reference to have the same order, the controller to be static, or state measurements to be available.

\textbf{Paper organization}:
Section~\ref{sec:preliminaries} provides preliminary results from behavioral systems theory.
Section~\ref{sec:problem_formulation} formalizes the data-driven controller implementability problem.
Section~\ref{sec:main_results} contains the main results of the paper, including necessary and sufficient implementability conditions which rely only on raw data and a formula for the direct data-driven construction of controllers that implement any plant-compatible reference behavior.
Section~\ref{sec:conclusion} provides a summary and an outlook to future research directions. The proofs of our main results are deferred to the appendix. 

\textbf{Notation}:  
The set of positive integers is denoted by $\N$. The set of real numbers is denoted by $\R$. For ${T\in\N},$ the set of integers $\{1, 2, \dots , T\}$ is denoted by $\mathbf{T}$.  The   image,  kernel, and Moore-Penrose inverse  of the matrix ${M \in \R^{p \times m}}$ are denoted by $\Image M$, $\ker M$, and $M^\dagger,$ respectively.  The collection of all maps from $X$ to $Y$ is denoted by $(Y)^{X}$.  The inverse image of the set $Y$ under $f$ is denoted by $f^{-1}(Y)$.

\section{Preliminary results}\label{sec:preliminaries}

This section recalls key notions and results from behavioral systems theory~\cite{willems1986timeI}, with a focus on discrete-time LTI systems.

\subsection{Time series and Hankel matrices}

We use the terms \textit{time series} and \textit{trajectory} interchangeably. The set of time series ${w=(w(1),\ldots,w(T))}$ of length ${T\in\N}$, with ${w(t)\in\R^{q}}$ for ${t\in\mathbf{T}}$, is defined as  ${(\R^q)}^\mathbf{T}$. The set of infinite-length time series $w=(w(1), w(2), \ldots)$, with  ${w(t)\in\R^{q}}$ for ${t\in\N}$, is defined as ${(\R^q)}^\N$. 

\subsubsection{The cut operator}
Restricting time series over subintervals gives rise to the cut operator. Formally, given 
${w\in {(\R^q)}^\mathbf{T}}$ and ${L\in\mathbf{T}}$, the \textit{cut operator} is defined as
\beq \nn
w|_L = (w(1),\ldots,w(L)) \in {(\R^q)}^\mathbf{L} .
\eeq
For infinite-length time series, the definition holds verbatim with ${w\in{(\R^q)}^\N}$ and ${L\in\N}$. Applied to a set of time series ${\mathcal{W} \subseteq {(\R^q)}^\mathbf{T}}$ or ${\mathcal{W} \subseteq {(\R^q)}^\N}$, the cut operator acts on all time series, defining the \emph{restricted} set $\mathcal{W}|_{L}  =\{w|_L \,:\, w\in \mathcal{W} \}$. By a convenient abuse of notation, we identify the trajectory $w|_L$ with the corresponding vector ${(w(1),\ldots,w(L))\in\R^{qL}}$.

\subsubsection{The shift operator} 
Shifting elements of time series gives rise to the shift operator.  Formally, given 
${w\in{\R^{qT}}}$ and $\tau\in\mathbf{T}$, the \textit{shift operator} is defined as
\beq \nn
\sigma^{\tau-1}w = (w(\tau),\ldots,w(T)) \in {\R^{q(T-\tau+1)}} .
\eeq
For infinite-length time series, the shift operator is defined as ${w \mapsto \sigma^{\tau-1}w}$, with ${\sigma^{\tau-1}w(t) = w(t+\tau-1)}$, for any ${\tau \in \N}$. Applied to a set of time series ${\mathcal{W} \subseteq {\R^{qT}}}$ or ${\mathcal{W} \subseteq {(\R^q)}^\N}$,
the shift operator acts on all time series in the set giving rise to the \emph{shifted} set $\sigma^\tau\mathcal{W}  =\{\sigma^\tau w\,:\, w\in\mathcal{W} \}$.

\subsubsection{Hankel matrices}
%
The \textit{Hankel matrix} of depth ${L\in\mathbf{T}}$ associated with the time series ${w \in \R^{qT}}$ is defined as 
\beq \label{eq:Hankel} 
\! H_{L}(w) \! = \!
\scalebox{0.85}{$
\bma  \nn
\begin{array}{ccccc}
w(1) & w(2)  & \cdots &  w(T-L+1)   \\
w(2) & w(3)  & \cdots &   w(T-L+2)   \\
\vdots  & \vdots  & \ddots & \vdots  \\
w(L) & w(L+1)  & \cdots  & w(T)
\end{array}
\ema
$} . \! \!
\eeq

\subsection{Discrete-time LTI dynamical systems}

A \textit{dynamical system} (or, briefly, \textit{system}) is a triple $\Sigma=(\T,\W,\B),$ where $\T$ is the \textit{time set}, $\W$ is the \textit{signal space}, and $\B \subseteq (\W)^{\T}$ is the \textit{behavior} of the system. 
We exclusively focus on \textit{discrete-time} systems, with ${\T = \N}$ and ${\W = \R^q}$.

\subsubsection{Finite-dimensional LTI systems}
A system $\B$ is \textit{linear} if $\B$ is a linear subspace, \textit{time-invariant} if $\B$ is shift-invariant, \textit{i.e.}, ${\sigma^{\tau-1}(\B) \subseteq \B}$ for all ${\tau \in \N}$, and \textit{complete} if $\B$ is closed in the topology of pointwise convergence~\cite[Proposition 4]{willems1986timeI}.
The model class of all complete LTI systems is denoted by $\L^q$.  By a convenient abuse of notation, we write $\B \in \L^{q}$.

\subsubsection{Kernel representations}

Every finite-dimensional LTI system ${\B \in \L^{q}}$ admits a \textit{kernel representation} of the form
\beq \nn
\B = \ker R(\sigma) ,
\eeq
where the operator $R(\sigma)$ is defined by the polynomial matrix $R(z) = R_0 +R_1 z +\ldots+ R_{\ell}z^\ell,$ with $R_i\in\R^{p\times q}$ for $i\in\Bell,$ and the set $\ker R(\sigma)$ is defined as $\{w \,:\, R(\sigma)w = 0\}$. Without loss of generality, we assume that  $\ker R(\sigma)$  is a \emph{minimal}  kernel representation  of $\B$,  \textit{i.e.}, $p$ is as small as possible over all kernel representations of $\B$.

\subsubsection{Integer invariants of an LTI system}

The structure of  an LTI  system ${\B\in\L^q}$ is characterized by a set of integer invariants~\cite[Section 7]{willems1986timeI}, defined as
\begin{itemize}
\item the \textit{number of inputs} ${m(\B) = q-\text{row dim} R}$,  
\item the \textit{number of outputs} ${p(\B) = \text{row dim} R}$,  
\item the \textit{lag} ${\ell(\B) = \max_{i \in\mathbf{p}}\{\deg\text{row}_i R \}}$, and
\item the \textit{order} ${n(\B)=\sum_{i \in\mathbf{p}} \deg\text{row}_i R }$, 
\end{itemize}
where   
$\ker R(\sigma)$  is a minimal  kernel representation of $\B$, while   
${\text{row dim}R}$ and ${\deg\text{row}_i R}$ are the number of rows and the degree of the $i$-th row of $R(z)$, respectively. 
The integer invariants are intrinsic properties of a system, as they do not depend on its representation~\cite[Proposition X.3]{willems1991paradigms}. 

\subsubsection{Partitions}
Given a permutation matrix ${\Pi\in\R^{q\times q}}$ and an integer $0 < m < q$, the map
\beq \label{eq:partition} 
(u,y)= \Pi^{-1} w 
\eeq
defines a \textit{partition} of ${w\in\R^q}$ into the variables ${u\in\R^m}$ and ${y\in\R^{q-m}}$.  We write $w \sim (u,y)$ if~\eqref{eq:partition} holds for some  permutation matrix ${\Pi\in\R^{q\times q}}$ and integer ${0 < m < q}$.  
Any partition~\eqref{eq:partition} induces the natural  projections  ${ \pi_u:  w \mapsto u}$ and ${ \pi_y:  w \mapsto y}$.  
We call $(u,y)$ a \textit{partition of ${\B\in\L^{q}}$} if~\eqref{eq:partition} holds for all ${w\in\B}$.  

\subsubsection{State-space representations}
Every finite-dimensional LTI  system ${\B \in \L^{q}}$ can be described by the equations
\beq \label{eq:state-space}
\sigma x = Ax + Bu, \quad y=Cx+Du,
\eeq
and admits a (\textit{minimal}) \textit{input/state/output representation} 
\beq \label{eq:partition-ISO} 
\!  \! \B \!  = \!
\left\{
(u,y) \in (\R^{q})^\N \,:\, \exists \,x\in(\R^n)^{\N} \, \textup{s.t.}~\eqref{eq:state-space}~\text{holds}
\right\},
\eeq
where
$\scalebox{0.75}{$\bma\!
\begin{array}{cc}
A & B \\
C & D
\end{array}\!\ema $}
 \in \R^{(n+p)\times (n+m)}$ and  $m$, $n$, and $p$ are the number of inputs, the order, and the number of outputs of $\B$, respectively. 

\subsection{Data-driven representations of LTI systems}

The restricted behavior of a finite-dimensional,  discrete-time,  LTI system can be   represented as the image of a raw data matrix.  We summarize a version of this principle known as the \textit{fundamental lemma}~\cite{willems2005note}.

\begin{lemma}~\cite[Corollary 19]{markovsky2020identifiability} \label{lemma:fundamental_generalized}
Let ${\B \in \L^{q}}$ and ${w \in \B|_{T}}$. Assume ${\ell(\B)<L \le T}$.
Then  $\BL = \Image H_L(w)$  if and only if
\beq  \label{eq:generalized_persistency_of_excitation} 
\rank H_L(w)   =  m(\B)L+ n(\B).
\eeq
\end{lemma}

\noindent
The rank condition~\eqref{eq:generalized_persistency_of_excitation} is referred to as the  \emph{generalized persistency of excitation} condition~\cite{markovsky2020identifiability}. Thus, we call a trajectory ${w \in \B|_{T}}$ of a system ${\B \in \L^{q}}$ \emph{generalized persistently exciting (GPE) of order $L$} if~\eqref{eq:generalized_persistency_of_excitation} holds.
Different variations of this principle can be formulated under  a range of assumptions, see, \textit{e.g.}, the recent survey~\cite{markovsky2021behavioral} for an overview.

\section{Problem formulation} \label{sec:problem_formulation}

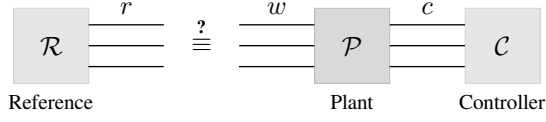
\begin{figure}[H]
\centering
\begin{tikzpicture} [scale=1, every node/.style={transform shape}]
\node at (12.5,-1) {${\mathcal{R}}$};
\filldraw[draw = black, opacity = 0.1]  (12,-0.5) rectangle (13,-1.5);
\node at (16.5,-1) {${\mathcal{P}}$};
\filldraw[draw = black, opacity = 0.15]  (16,-0.5) rectangle (17,-1.5);
\node at (18.5,-1) {${\mathcal{C}}$};
\filldraw[draw = black, opacity = 0.1]  (18,-0.5) rectangle (19,-1.5);
\draw [line width = .5 pt](12.9972,-1) -- (13.9972,-1);
\draw [line width = .5 pt](12.9952,-1.278) -- (13.9952,-1.278);
\draw [line width = .5 pt](12.9972,-0.722) -- (13.9972,-0.722); 
\draw [line width = .5 pt](14.9984,-1) -- (15.9984,-1);
\draw [line width = .5 pt](14.9964,-1.278) -- (15.9964,-1.278);
\draw [line width = .5 pt](14.9984,-0.722) -- (15.9984,-0.722);
\draw [line width = .5 pt](17,-1) -- (18,-1);
\draw [line width = .5 pt](16.998,-1.278) -- (17.998,-1.278);
\draw [line width = .5 pt](17,-0.722) -- (18,-0.722);
\node at (13.4984,-0.5) {${r}$};
\node at (15.5,-0.5) {${w}$};
\node at (17.5,-0.5) {${c}$};
\node at (14.5,-0.85) {$\stackrel{\textbf{?}}{\equiv}$}; 
\node at (12,-0.5) {};
\node at (12,-2) {};
\node at (19,-2) {};
\node at (19,-0.5) {};
\node at (12.5,-1.75) {\footnotesize{Reference}};
\node at (16.5,-1.75) {\footnotesize{Plant}};
\node at (18.5,-1.75) {\footnotesize{Controller}};
\end{tikzpicture}
\centering
\caption{Control in a behavioral setting.}
\label{fig:control}
\end{figure}

Consider a \textit{plant behavior} ${\Plant \in \L^{q+k}}$, a \textit{reference behavior} ${\Reference \in \L^{q}}$, and a \textit{controller behavior} ${\Controller \in \L^{k}}$, as shown in Fig.~\ref{fig:control}.

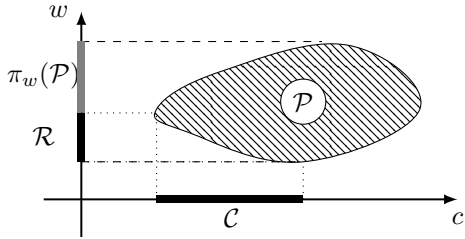
\begin{figure}[h!]
\centering
\usetikzlibrary{arrows}
\usetikzlibrary{patterns}
\begin{tikzpicture} [scale=1, every node/.style={transform shape}]
\draw [latex-, line width = .75 pt](14,3) -- (14,0); 
\draw [latex-, line width = .75 pt](19,0.5) -- (13.5,0.5); 
\node at (13.7,3) {${w}$};
\node at (19,0.25) {${c}$};
\node at (16.95,1.8) {${\mathcal{P}}$};
\filldraw[pattern=north west lines, even odd rule]  plot[smooth cycle, tension=.7] coordinates {(15,1.7) (15.4,1.4) (16.95,1) (18.5,1.7) (17.55,2.55) (15.85,2.2)} (16.95,1.8) circle (0.3);
\node at (13.5,2.15) {${\pi_w(\mathcal{P})}$};
\draw[dashed] (14,2.6) -- (17.2,2.6);
\draw [dashed](14,1) -- (16.95,1);
\fill [gray](13.95,1) -- (14.05,1) -- (14.05,2.6) -- (13.95,2.6)  -- cycle;
\node at (13.5,1.35) {${\mathcal{R}}$};
\fill [black](13.95,1) -- (14.05,1) -- (14.05,1.65) -- (13.95,1.65)  -- cycle; 
\node at (16,0.25) {${\mathcal{C}}$};
\fill [black](15,0.55) -- (15,0.45) -- (16.95,0.45) -- (16.95,0.55)  -- cycle;
\draw [dotted](15,0.5) -- (15,1.65);
\draw [dotted](16.95,0.5) -- (16.95,1);
\draw [dotted](14,1) -- (16.95,1);
\draw [dotted](14,1.65) -- (15,1.65);
\node at (13.5,3) {};
\node at (13.5,0) {};
\node at (19,0) {};
\node at (19,3) {};
\end{tikzpicture}
\centering
\caption{Pictorial illustration of behaviors ${\Plant}$, ${\Reference}$, ${\Controller}$, and $\pi_w(\Plant)$.}
\label{fig:behaviors}
\end{figure}

\noindent 
Following~\cite{willems2002synthesis}, we fix a partition of the variables of the plant behavior $\Plant$, which induces the natural projections 
\beq 
\pi_{w}:  (w,c) \mapsto w, \quad  \pi_{c}:  (w,c) \mapsto c  ,
\eeq
where $w$ are the \emph{to-be-controlled} variables and $c$ are the \emph{control} variables, respectively. 
The controller behavior $\Controller$ is interconnected to the plant behavior $\Plant$ via \textit{variable sharing}.
Formally, the \textit{interconnection of $\Plant$ and $\Controller$ via the shared variable $c$} is defined as
\beq \label{eq:composition}
\PC = 
\left\{ 
    (w,c) \in (\R^{q+k})^\N 
        \,:\, \, 
    c \in \mathcal{C}, \ (w,c) \in \mathcal{P}
\right\} . 
\eeq
Similarly, we define the \textit{hidden behavior} $\Hidden $ of the plant behavior $\Plant$ as
\beq \label{eq:hidden_behavior} 
\Hidden = 
\left\{ 
	w \in (\R^{q+k})^\N  \,:\,  \,  (w,0) \in \Plant
\right\} .
\eeq  
We refer to $\pi_w(\mathcal{P} \lVert_{c} \mathcal{C})$ and $\pi_w(\Plant)$ as the 
\textit{controlled plant behavior} and \textit{uncontrolled plant behavior}, respectively.  
Fig.~\ref{fig:behaviors} offers a pictorial illustration of the aforementioned behaviors. 
\vspace{-0.15cm}
\vspace{-0.15cm}


A controller ${\Controller \in \L^{k}}$ is said to \textit{implement} ${\Reference \in \L^{q} }$ if ${\pi_w(\PC)  = \Reference}$~\cite{willems2002synthesis}. In other words, a controller behavior implements a given reference behavior if the resulting controlled plant behavior obtained from interconnecting the plant with the controller coincides with the reference behavior.  Consequently, a behavior ${\Reference \in \L^{q} }$ is said to be \textit{implementable} if there exists a controller which implements $\Reference$. 
 
  
\begin{problem}[Data-driven  controller  implementability problem]
 Consider a plant behavior  ${\Plant \in \L^{q+k}}$ and a reference behavior  ${\Reference \in \L^{q}}$. 
Given trajectories of length ${T\in\N}$ of the plant behavior ${(w,c) \in \Plant|_{T}}$ and of the reference behavior ${r \in \Reference|_{T}}$, the \textit{data-driven controller implementability problem} is to find, if possible, a  controller   ${\Controller \in \L^k}$ which implements  $\Reference$.
\end{problem}

\noindent
The data-driven implementability problem is \textit{solvable} if ${\Reference }$ is implementable, in which case any controller ${\Controller}$ implementing ${\Reference }$ is a \textit{solution} of the problem.  


\section{Main results}\label{sec:main_results}

This section contains the main results of the paper and is logically divided in two parts. First, we provide necessary and sufficient conditions for implementablity of a given reference behavior which only rely on measured data. Second, we present a data-driven strategy to obtain controllers  for any   given implementable reference behavior.

\subsection{Data-driven implementability conditions}

The data-driven  controller implementability  problem is closely related to  the controller implementability problem~\cite{willems2002synthesis}, which seeks to determine all implementable reference behaviors ${\Reference \in \L^q}$ for a given plant $\Plant\in\L^{q+k}$.  
The  problem  has been first studied in a continuous-time setting in~\cite{willems2002synthesis}. An elegant solution is provided by the following result.

\begin{theorem}[Infinite-horizon implementability conditions]\cite[Theorem 1]{willems2002synthesis} \label{thm:controller_implementability}
 Consider a plant behavior  ${\Plant \in \L^{q+k}}$ and a reference behavior  ${\Reference \in \L^{q}}$.  Then ${\Reference}$ is implementable if and only if
\beq \label{eq:implementability} 
\Hidden \subseteq \Reference \subseteq \pi_w(\Plant).
\eeq
\end{theorem} 
\noindent  
Theorem 1 provides a powerful necessary and sufficient condition for the existence of controllers implementing a given reference behavior.  
However, verifying the implementability condition~\eqref{eq:implementability} 
 may be challenging  in practice because it requires full knowledge of both the hidden behavior and the uncontrolled plant behavior;  this is especially true if only measured data of the plant and reference behaviors are  available.


We now present a simple, but important extension of the controller implementability theorem, which 
provides necessary and sufficient condition for a reference behavior to be implementable while only requiring knowledge of the hidden behavior and the uncontrolled plant behavior over a \textit{finite} time horizon. 

\begin{theorem}[Finite-horizon implementability conditions] \label{thm:implementability-finite-horizon}
Consider a plant behavior ${\Plant \in \L^{q+k}}$   and a reference behavior  ${\Reference \in \L^{q}}$.  Suppose ${L>\max\{\ell(\Plant),\ell(\Reference),\ell(\pi_w(\Plant))\}}$. Then ${\Reference}$ is implementable if and only if
 \beq \label{eq:implementability-finite-horizon} 
\HiddenL \subseteq \ReferenceL \subseteq \pi_w(\Plant)|_L.
\eeq
\end{theorem}

\noindent  
Theorem~\ref{thm:implementability-finite-horizon} offers an alternative non-parametric necessary and sufficient condition for the existence of controllers that can implement a given reference behavior.  Similar to Theorem~\ref{thm:controller_implementability}, Theorem~\ref{thm:implementability-finite-horizon} establishes implementability conditions that do not rely on a specific representation. However, unlike Theorem~\ref{thm:controller_implementability}, the subspace inclusions~\eqref{eq:implementability-finite-horizon} only need information about finite-horizon behaviors, whereas subspace inclusions~\eqref{eq:implementability} require knowledge of the complete (infinite-dimensional) behaviors.


An important consequence of Theorem~\ref{thm:implementability-finite-horizon} is that the subspace inclusions~\eqref{eq:implementability-finite-horizon} can be translated into implementability criteria which can be verified directly from data.  
In particular, the following result provides general necessary and sufficient conditions for the implementability of a given reference behavior using data. 

\begin{corollary}[Data-driven implementability conditions]\label{cor:behaviors-from-data}
Consider a plant behavior  ${\Plant \in \L^{q+k}}$ and a reference behavior  ${\Reference \in \L^{q}}$. Suppose ${L>\max\{\ell(\Plant),\ell(\Reference),\ell(\pi_w(\Plant))\}}$. Let  ${(w,c)\in\Plant|_T}$ and ${r\in\Reference|_T}$ be GPE of order $L$.  Define 
\begin{align} 
N &=  H_L(w)\left(I-H_L(c)^{\dagger}H_L(c)\right)\label{eq:Bhidden-data} \\
R &=  H_L(r) \label{eq:Br-data} \\
 P_w  &= H_L(w)\label{eq:Buncontrolled-data}.
\end{align}
Then 
${\HiddenL = \Image N}$, 
${\ReferenceL = \Image R}$, 
${\pi_w(\Plant)|_L = \Image  P_w }$.
Consequently, the reference behavior $\Reference$ is implementable if and only if the system of linear equations
\beq \label{eq:EM-solvability-data} 
N = R \Phi , \quad R   = P_w \Psi ,  
\eeq 
in the unknown matrices $\Phi$ and $\Psi$  admits a solution. 
\end{corollary}
\noindent
Corollary~\ref{cor:behaviors-from-data} establishes necessary and sufficient conditions for testing the implementability of a given reference behavior directly from data. This, in turn, provides a necessary and sufficient  condition for the solvability of the data-driven controller implementability problem.

\subsection{Data-driven canonical controller representation}

Theorem~\ref{thm:implementability-finite-horizon} and Corollary~\ref{cor:behaviors-from-data} provide conditions under which a reference behavior is implementable, but do not provide expressions for a controller which implements the reference behavior.
 We first recall an expression for a controller $\Controller \in \L^{k}$ which implements a given implementable reference behavior $\Reference \in \L^{q}$ and, subsequently,  obtain an expression for such controller which relies only on data.  

\begin{theorem}[Canonical controller]\label{thm:canonical-controller}~\cite[Theorem 2.1]{van2003achievable}
Consider a plant behavior  ${\Plant \in \L^{q+k}}$ and a reference behavior  ${\Reference \in \L^{q}}$. Assume $\Reference$ is implementable. Then $\Reference$ is implemented by the controller
\beq \label{eq:canonical_controller}
\Controller = \pi_{c}(\PR).
\eeq 
\end{theorem}

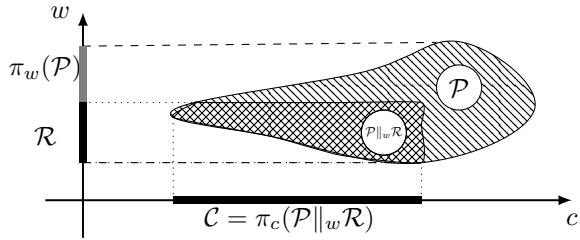
\begin{figure}[h!]
\centering
\begin{tikzpicture} [black, scale=1, every node/.style={transform shape}]
\draw [latex-, line width = .75 pt](14,3) -- (14,0); 
\draw [latex-, line width = .75 pt](20.5,0.5) -- (13.5,0.5); 
\node at (13.7,3) {${w}$};
\node at (20.5,0.25) {${c}$};
\node at (19,2) {${\mathcal{P}}$};
\filldraw[pattern=north west lines, even odd rule]  plot[smooth cycle, tension=.7] coordinates {(15.2,1.7) (16.5,1.35) (18.5,1) (20,1.7) (19.05,2.6) (17.7,2.2)} (19,2) circle (0.3) (18,1.4) circle (0.3);
\node at (13.5,2.25) {${\pi_w(\mathcal{P})}$};
\draw[dashed] (14,2.55) -- (18.7,2.6);
\draw [dashed](14,1) -- (18.5,1);
\fill [gray](13.95,1) -- (14.05,1) -- (14.05,2.55) -- (13.95,2.55)  -- cycle;
\node at (13.5,1.35) {${\mathcal{R}}$};
\fill [black](13.95,1) -- (14.05,1) -- (14.05,1.8) -- (13.95,1.8)  -- cycle; 
\node at (16.75,0.25) {${\mathcal{C} = \pi_{c}(\PR)}$};
\fill [black](15.2,0.55) -- (15.2,0.45) -- (18.5,0.45) -- (18.5,0.55)  -- cycle;
\draw [dotted](15.2,0.5) -- (15.2,1.65);
\draw [dotted](18.5,0.5) -- (18.5,1);
\draw [dotted](14,1) -- (18.5,1);
\draw [dotted](14,1.8) -- (15.5,1.8);
\node at (18,1.4) {\scalebox{0.5}{${\mathcal{P} \lVert_w \mathcal{R}}$}};
\filldraw[pattern=north east lines, even odd rule]  plot[smooth cycle, tension=.3 ] coordinates {(17,1.8) (16,1.8)  (15.5,1.8)  (15.5,1.8)  (15.2,1.6)  (16.5,1.35) (17.42,1.1)(18.5,1) (18.5,1.5) (18.5,1.8) (17.5,1.8)}  (18,1.4) circle (0.3);
\node at (13.5,3) {};
\node at (13.5,0) {};
\node at (20.5,0) {};
\node at (20.5,3) {};
\end{tikzpicture}
\centering
\caption{Pictorial illustration of the canonical controller.}
\label{fig:canonical_controller}
\end{figure}

\noindent 
Theorem~\ref{thm:canonical-controller} provides a universal formula which defines the behavior of a controller which implements any implementable reference behavior. Consequently, the controller~\eqref{eq:canonical_controller} is referred to as the \textit{canonical controller}. 
Fig.~\ref{fig:canonical_controller} offers a pictorial illustration of the canonical controller behavior. 

\vspace{-0.15cm}
\vspace{-0.15cm}

The concept of canonical controller has been implicitly defined in the seminal paper~\cite{willems2002synthesis} and extended to general systems, e.g., in~\cite{van2003achievable} and~\cite{julius2005canonical}. The canonical controller is appealing due to its simple construction and its representation-free formalization of the internal model principle~\cite{francis1976internal}. 

\begin{remark}[Canonical controller and well-posedness]
The canonical controller is such that the interconnection   of $\Plant$ and $\Controller$  is  well-posed, \textit{i.e.}, $\PC\neq \emptyset$. If $\Reference$ is implementable, then 
$$\PC 
\stackrel{\eqref{eq:canonical_controller}}{=} \Plant \cap  (\PR) 
\stackrel{\eqref{eq:implementability} }{=} \PR \neq \emptyset , $$ 
where the last inequality follows from the implementability assumption on $\Reference$.
\end{remark}

Next, we show that this concept also allows us to define the restricted behavior of a controller ${\Controller \in \L^{k}}$ which implements a desired reference behavior ${\Reference \in \L^{q}}$ using only measured data. 
For ${L\in\N}$, we define 
the matrix representations ${\Pi_{w} \in \R^{qL \times (q+k)L}}$ and  ${\Pi_{c} \in \R^{kL \times (q+k)L}}$ of the projections $\pi_w$ and $\pi_{c}$ over the time horizon ${[1,L]}$ as
\begin{align*}
\Pi_{w}  &= 
\blockdiag\left(
\begin{bmatrix}
I & 0 
\end{bmatrix},
\dots,
\begin{bmatrix}
I & 0 
\end{bmatrix}
 \right) , \\
\Pi_{c}  &= 
\blockdiag\left(
\begin{bmatrix}
0 & I
\end{bmatrix},
\dots,
\begin{bmatrix}
0 & I
\end{bmatrix}
 \right) .
\end{align*}

\begin{corollary}[Data-driven canonical controller representation] \label{cor:controller-behavior-from-data} 
Consider a plant behavior  ${\Plant \in \L^{q+k}}$ and a reference behavior  ${\Reference \in \L^{q}}$.  
Assume $\Reference$ is implementable and let $\Controller$ be the canonical controller~\eqref{eq:canonical_controller}. Let ${L>\max\{\ell(\Plant),\ell(\Reference),\ell(\pi_w(\Plant))\}}$. 
Let ${(w,c)\in\Plant|_T}$ and ${r\in\Reference|_T}$ be GPE of order $L$. Define 
\beq \nn
P \sim  H_L((w,c)) , \quad  R= H_L(r) ,
\eeq
and
\beq \nn
P_p = PP^{\dagger} ,  \quad 
P_r \sim 
\bma
\begin{array}{cc}
RR^\dagger & 0\\
0 & I_{kL}
\end{array}
\ema ,
\eeq
where $\sim$ denotes similarity under a coordinates permutation.
Then
\beq \label{eq:Bcontroller-data} 
\ControllerL =  \Image \Pi_{c} P_r\left(P_r+P_p\right)^\dagger P_p .
\eeq
\end{corollary}
 
Corollary~\ref{cor:controller-behavior-from-data} provides a data-based description of the finite-horizon behavior of the canonical controller. This formula serves a dual purpose: it can be used to identify a controller from measured data of the reference and the plant, or for direct control purposes by generating finite-length trajectories of the canonical controller.  Note that longer trajectories for specific control requirements may be also generated using the lemma on weaving trajectories~\cite[Lemma 8.21]{markovsky2006exact}.

\begin{remark}[Persistency of excitation of the data]
Corollaries~\ref{cor:behaviors-from-data} and~\ref{cor:controller-behavior-from-data} rely on the assumption that ${(w,c)\in\Plant|_T}$ and ${r\in\Reference|_T}$ are GPE of order $L$. In order to check such assumption from data, upper bounds on $n(\Plant)$, $m(\Plant)$, $n(\Reference)$, and $m(\Reference)$ are needed (see the rank condition~\eqref{eq:generalized_persistency_of_excitation}).  Alternatively, the rank condition~\eqref{eq:generalized_persistency_of_excitation} can be guaranteed to hold for controllable systems if a certain rank condition on the inputs hold~\cite{willems2005note}.
\end{remark}


\begin{remark}[Alternative matrix representations]
Corollaries~\ref{cor:behaviors-from-data} and~\ref{cor:controller-behavior-from-data} can be also expressed using alternative data-driven representations of the restricted behaviors $\HiddenL,$  $\ReferenceL$, and $\pi_w(\Plant)|_L$, \textit{e.g.}, using  Page matrices~\cite{coulson2021distributionally} or mosaic-Hankel matrices~\cite{van2020willems,markovsky2022identifiability}. 
\end{remark}

\begin{remark}[Connections to exact model matching]
The controller implementability problem is closely related to the \textit{exact model matching}~\cite{wolovich1972use,wang1972exact}, where the goal is to design a state feedback law for an LTI system to match a reference transfer function.
The problem is well-studied in a model-based context~\cite{wolovich1972use,wang1972exact,wang1972solution} and it typically reduces to solving a set of linear algebraic equations~\cite{wang1972exact,wang1972solution}. 
The recent paper~\cite{breschi2021direct} presents analogous findings in a data-driven context. It can be shown that the data-driven implementability condition~\eqref{eq:EM-solvability-data} generalizes the results obtained in~\cite{breschi2021direct}.
\end{remark}

\begin{remark}[Reference behaviors that are not implementable]
When the reference behavior $\Reference$ does not satisfy the implementability conditions~\eqref{eq:implementability} or, equivalently,~\eqref{eq:implementability-finite-horizon} for ${L>\max\{\ell(\Plant),\ell(\Reference),\ell(\pi_w(\Plant))\}}$, one option is to adjust $\Reference$~\cite[Remark 2.6]{van2003achievable}. This involves excluding $w$ values without corresponding $c$ values in $\Plant$ and including $w$ values from $\Reference$ that match with $c$ values in $\Plant$, creating a new implementable reference $\Reference^{\prime}$. Alternatively, one may  search for implementable controlled behaviors such that~\eqref{eq:implementability-finite-horizon} holds, while minimizing the distance from the original reference behavior by exploiting the (Grassmannian) geometry of finite-horizon LTI behaviors~\cite{padoan2019behavioral}.
\end{remark}

\section{Conclusion}\label{sec:conclusion}

We have studied the controller implementability problem from the lens of data-driven control, providing necessary and sufficient implementability conditions which can rely solely on raw data. Furthermore, we have addressed the problem of constructing controllers directly from data. By employing the notion of canonical controller, we have presented a formula for generating controllers which implement any plant-compatible reference behaviors in a data-driven fashion.  Future research should address noisy scenarios and study an approximate version of the controller implementability problem.

\appendix 
\subsection{Proofs}

\subsubsection{Preliminary results}
The proofs of our main results rely on several preliminary results about the interplay between coordinate projections, the cut operator, restricted LTI behaviors, and orthogonal projections onto intersections of subspaces.



\begin{lemma}[Preimages under coordinate projections] \label{lemma:preimage_of_behaviors_under_coordinate_projections}
Let ${\B\in \L^{q}}$. Assume $(w,c)$ is a partition of $(\R^{q+k})^\N$, with ${w \in \R^{q}}$ and ${c \in \R^{k}}$. Then
$\pi_{w}^{-1}(\B) = \B \times (\R^k)^{\N} .$
\end{lemma}

\begin{proof}
By definition, we have 
\begin{align*}
\pi_{w}^{-1}(\B) 
= \left\{ (w,c) \in (\R^{q+k})^{\N}  \,:\, w \in \B \right\} = \B \times (\R^k)^{\N} .
\end{align*}
\end{proof}

\begin{lemma}[Coordinate projections and LTI behaviors] \label{lemma:coordinate_projection_and_LTI_behaviors}
Let ${\B \in \L^{q+k}}$. Assume $(w,c)$ is a partition of $\B$, with ${w \in \R^{q}}$ and ${c \in \R^{k}}$. Then $\pi_{w}(\B)|_{L} = \Pi_{w}(\BL)$ for all ${L\in\N}$.
\end{lemma}

\begin{proof}
$(\subseteq)$. We first show ${\pi_{w}(\B)|_{L} \subseteq \Pi_{w}(\BL)}.$ 
Let ${(\tilde{w}(1),\dots,\tilde{w}(L))\in\pi_{w}(\B)|_L}$. Then there is ${w\in \pi_w(\B)}$ such that 
$w|_L=(\tilde{w}(1),\dots,\tilde{w}(L)).$
Furthermore, there exists $c$ such that ${(w,c)\in\B}$ and, hence, 
$$(w,c)|_L=(\tilde{w}(1),c(1),\dots,\tilde{w}(L),c(L))\in\BL.$$ 
This implies $(\tilde{w}(1),\dots,\tilde{w}(L))\in\Pi_{w}(\BL),$
and, hence, ${\pi_{w}(\B)|_L \subseteq \Pi_{w}(\B|_L)}$.

\noindent$(\supseteq)$. Next, we show 
${\pi_{w}(\B)|_{L} \supseteq \Pi_{w}(\BL)}.$ Let 
$${(\tilde{w}(1),\dots,\tilde{w}(L))\in\Pi_w(\BL)}.$$ 
Then there exists $(\tilde{c}(1),\dots,\tilde{c}(L))$ such that 
$$(\tilde{w}(1),\tilde{c}(1),\dots,\tilde{w}(L),\tilde{c}(L))\in\BL.$$
Thus, there exists ${(w,c)\in\B}$ such that 
$$(w,c)|_L=(\tilde{w}(1),\tilde{c}(1),\dots,\tilde{w}(L),\tilde{c}(L)).$$
Then ${w\in \pi_w(\B)}$ and, hence, ${w|_L\in\pi_{w}(\B)|_{L}}$. This  implies 
$${(\tilde{w}(1),\dots,\tilde{w}(L))\in\pi_{w}(\B)|_{L}},$$ 
which proves ${\Pi_{w}(\BL) \subseteq \pi_{w}(\B)|_{L}}$ and, hence, the claim.
\end{proof}

\begin{lemma}[Intersection of restricted LTI behaviors] \cite[Proposition 16]{yan2023distributed} \label{lemma:intersection_of_restricted_behaviors}
Let ${\B\in\L^{q}}$ and ${\bar{\B}\in\L^{q}}$. Then
\beq \nn
(\B \cap \bar{\B})|_L \subseteq \B|_L \cap \bar{\B}|_L.
\eeq
for all ${L\in\N}.$ Furthermore, if $L>\max\{\ell(\B), \ell(\bar{\B})\}$, then
\beq \nn
(\B \cap \bar{\B})|_L = \B|_L \cap \bar{\B}|_L.
\eeq
\end{lemma}

\begin{lemma}[Cartesian product of restricted LTI behaviors] \cite[Proposition 19]{yan2023distributed} \label{lemma:cartesian_product_of_LTI_behaviors}
Let ${\B \in \L^{q}}$ and ${\bar{\B} \in \L^{k} }$. Then for all ${L\in\N}$,
$(\B \times \bar{\B})|_L = \B|_L \times \bar{\B}|_L$.
\end{lemma}

\begin{lemma}[Inclusion between restricted LTI behaviors] \label{lemma:behavior_subsets} 
Let ${\B\in \L^{q}}$ and ${\bar{\B}\in \L^{q}}$. Then ${\bar{\B}|_L\subseteq \B|_L}$ implies ${\bar{\B}\subseteq \B}$ for $L>\max\{\ell(\bar{\B}),\ell(\B)\}$.
\end{lemma}
\begin{proof}
We have that ${\bar{\B}|_L\subseteq \B|_L}$ if and only if ${\bar{\B}|_L\cap \B|_L}=\bar{\B}|_L$. By Lemma~\ref{lemma:intersection_of_restricted_behaviors}, we have  $${\bar{\B}|_L\cap \B|_L}=(\B \cap \bar{\B})|_L.$$ Then $\bar{\B}|_L=(\B \cap \bar{\B})|_L$. By~\cite[Corollary 14]{markovsky2020identifiability}, $\bar{\B}=\B \cap \bar{\B}$ which implies that ${\bar{\B}\subseteq \B}$, proving the claim.
\end{proof}

\begin{lemma}[Projectors on intersection of subspaces] \cite[p.2]{ben2015projectors} \label{lemma:projectors_on_intersection_of_subspaces}
Let $\mathcal{V}$ and $\mathcal{W}$ be subspaces of $\R^{n}$ and let $P_\mathcal{V}$ and $P_\mathcal{W}$ be the orthogonal projectors on $\mathcal{V}$ and $\mathcal{W}$, respectively. Then the orthogonal projector on the intersection of $\mathcal{V}$ and $\mathcal{W}$ is
\beq \nn
P_{\mathcal{V} \cap \mathcal{W}} = 2 P_\mathcal{V} (P_\mathcal{V}+P_\mathcal{W})^{\dagger} P_\mathcal{W} .
\eeq
\end{lemma}

\subsubsection{Proof of Theorem~\ref{thm:implementability-finite-horizon}}
We prove the claim by showing that~\eqref{eq:implementability} is equivalent to~\eqref{eq:implementability-finite-horizon} for ${L>\max\{\ell(\Plant),\ell(\Reference),\ell(\pi_w(\Plant))\}}$. 

\noindent\eqref{eq:implementability}$\implies$\eqref{eq:implementability-finite-horizon}: This holds by definition of the cut operator.

\noindent\eqref{eq:implementability-finite-horizon}$\implies$\eqref{eq:implementability}: We first show that $\HiddenL\subseteq\ReferenceL$ implies $\Hidden\subseteq\Reference.$ First, note that $L>\max\{\ell(\Hidden),\ell(\Reference)\}$. Indeed, by assumption, $L>\ell(\Reference)$. Furthermore, $L> \ell(\Hidden)$. Indeed, let $$\Plant=\ker\begin{bmatrix} R_w(\sigma) & R_c(\sigma)\end{bmatrix}$$ be a minimal kernel representation for $\Plant$. Then $\Hidden=\ker R_w(\sigma)$. Thus, $\ell(\Plant)\geq \ell(\Hidden)$ and, hence, $L>\max\{\ell(\Hidden),\ell(\Reference)\}$.  By Lemma~\ref{lemma:behavior_subsets}, we conclude that $\Hidden\subseteq\Reference$.
It can be shown that $\ReferenceL\subseteq\pi_w(\Plant)|_L$ implies $\Reference\subseteq\pi_w(\Plant)$ using similar arguments. \qed

\subsubsection{Proof of Corollary~\ref{cor:behaviors-from-data}}
Let $\bar{w}\in\pi_w(\Plant)|_L$. Since ${(w,c)\in\Plant|_T}$ and ${r\in\Reference|_T}$ are GPE of order $L$, there exists $g$ and $\bar{c}$ such that 
\[
\begin{bmatrix}
H_L(w)\\
H_L(c)
\end{bmatrix}g
=
\begin{bmatrix}
\bar{w}\\
\bar{c}
\end{bmatrix}.
\]
Thus $\bar{w}\in \Image H_L(w)$. Now let $\bar{w}\in \Image H_L(w)$. Then there exists $g$ such that $H_L(w)g=\bar{w}$. Thus there exists $\bar{c}\in \Image H_L(c)$ such that $(\bar{w},\bar{c})\in\PlantL$. Hence, $\bar{w}\in\pi_w(\Plant)|_L$, so that $\pi_w(\Plant)|_L=\Image H_L(w)$.

Now let $\bar{w}\in \HiddenL$. Then there exists $g$ such that
\[
\begin{bmatrix}
H_L(w)\\
H_L(c)
\end{bmatrix}g
=
\begin{bmatrix}
\bar{w}\\
0
\end{bmatrix}.
\]
Thus, ${g\in \ker H_L(c)=\Image (I-H_L(c)^\dagger H_L(c))}$. This, in turn, implies $$\bar{w}\in \Image H_L(w)(I-H_L(c)^\dagger H_L(c)).$$
Now let $\bar{w}\in \Image H_L(w)(I-H_L(c)^\dagger H_L(c))$. Then there exists $g$ such that $\bar{w}= H_L(w)(I-H_L(c)^\dagger H_L(c))g$. Thus,
\[
\begin{bmatrix}
\bar{w}\\
0
\end{bmatrix}=
\begin{bmatrix}
H_L(w)(I-H_L(c)^\dagger H_L(c))\\
0
\end{bmatrix}g\\
=\begin{bmatrix}
H_L(w)\\
H_L(c)
\end{bmatrix}\bar{g},
\]
with $\bar{g}=(I-H_L(c)^\dagger H_L(c))g$.
Hence, $(\bar{w},0)\in\PlantL$, so $\bar{w}\in\HiddenL$, showing that $$\HiddenL=\Image (I-H_L(c)^\dagger H_L(c)).$$ The second claim now follows directly from Theorem~\ref{thm:implementability-finite-horizon} and  the fact that the subspace inclusions
$$\Image N \subseteq \Image R  \subseteq \Image P_w $$
hold if and only if the system of linear equations~\eqref{eq:EM-solvability-data} admits a solution.\qed

\subsubsection{Proof of Corollary~\ref{cor:controller-behavior-from-data}}

By assumption, $\Reference$ is implementable and, hence, the canonical controller~\eqref{eq:canonical_controller} is well-defined.
By applying the cut operator to the definition of the canonical controller~\eqref{eq:canonical_controller}, we obtain
\[
\Controller|_L=\pi_{c}( \pi_{w}^{-1}(\Reference) \cap \Plant)|_L.
\]
By Lemma~\ref{lemma:coordinate_projection_and_LTI_behaviors},  we obtain
$
\Controller|_L=\Pi_{c}\left( (\pi_{w}^{-1}(\Reference) \cap \Plant)|_L \right).
$
Using Lemma~\ref{lemma:intersection_of_restricted_behaviors} and  ${{L}>\max\{\ell(\Plant),\ell(\Reference),\ell(\pi_w(\Plant))\}}$,  we obtain
\[
\Controller|_L=\Pi_{c}\left(\pi_{w}^{-1}(\Reference)|_L \cap \Plant|_L \right).
\]
By Lemma~\ref{lemma:preimage_of_behaviors_under_coordinate_projections} and
Lemma~\ref{lemma:cartesian_product_of_LTI_behaviors}, we can write the above as
\[
\Controller|_L=\Pi_{c}\left((\Reference|_L \times \R^{kL}) \cap \Plant|_L \right).
\]
By Lemma~\ref{lemma:projectors_on_intersection_of_subspaces}, the fact that ${(w,w_c)\in\Plant|_T}$ and ${w_r\in\Reference|_T}$ are GPE of order $L$, 
and the definition of the projectors  $P_p$ and $P_r$, we obtain $\Controller|_L=\Pi_{c} \Image  P_r\left(P_r+P_p\right)^\dagger P_p.$
Finally, since $\Pi_{c}$ is surjective,
\[
\Controller|_L= \Image \Pi_{c} P_r\left(P_r+P_p\right)^\dagger P_p. 
\]
This proves the result.\qed

\bibliographystyle{IEEEtran}
\bibliography{refs}

\newcommand{\TAC}{\textit{{IEEE} Trans. Autom.
  Control}}\newcommand{\TCST}{\textit{{IEEE} Trans. Syst.
  Tech.}}\newcommand{\TIT}{\textit{{IEEE} Trans. Inform.
  Theory}}\newcommand{\TSP}{\textit{{IEEE} Trans. Sign.
  Proc.}}\newcommand{\SCL}{\textit{Syst. Control
  Lett.}}\newcommand{\IJC}{\textit{Int. J.
  Control}}\newcommand{\EJC}{\textit{Eur. J.
  Control}}\newcommand{\ACC}{\textit{Proc. Amer. Control
  Conf.}}\newcommand{\ECC}{\textit{Proc. Eur. Control
  Conf.}}\newcommand{\CDC}[1]{\textit{Proc. {#1} Conf. Decision
  Control}}\newcommand{\CDCs}[1]{\textit{{#1} Conf. Decision
  Control}}\newcommand{\IFAC}[1]{\textit{Proc. {#1} IFAC World
  Congr.}}\newcommand{\SIAM}{\textit{SIAM J. Control
  Optim.}}\newcommand{\MCSS}{\textit{Math. Control, Sign.
  Syst.}}\newcommand{\NOLCOS}[1]{\textit{Proc. {#1} IFAC Symp. Nonlinear
  Control Syst.}}\newcommand{\MTNS}[1]{\textit{ Proc. {#1} Math. Symp. Netw.
  Syst.}}
\begin{thebibliography}{10}
\providecommand{\url}[1]{#1}
\csname url@samestyle\endcsname
\providecommand{\newblock}{\relax}
\providecommand{\bibinfo}[2]{#2}
\providecommand{\BIBentrySTDinterwordspacing}{\spaceskip=0pt\relax}
\providecommand{\BIBentryALTinterwordstretchfactor}{4}
\providecommand{\BIBentryALTinterwordspacing}{\spaceskip=\fontdimen2\font plus
\BIBentryALTinterwordstretchfactor\fontdimen3\font minus
  \fontdimen4\font\relax}
\providecommand{\BIBforeignlanguage}[2]{{%
\expandafter\ifx\csname l@#1\endcsname\relax
\typeout{** WARNING: IEEEtran.bst: No hyphenation pattern has been}%
\typeout{** loaded for the language `#1'. Using the pattern for}%
\typeout{** the default language instead.}%
\else
\language=\csname l@#1\endcsname
\fi
#2}}
\providecommand{\BIBdecl}{\relax}
\BIBdecl

\bibitem{willems1997interconnections}
J.~C. Willems, ``On interconnections, control, and feedback,'' \emph{IEEE
  Transactions on Automatic control}, vol.~42, no.~3, pp. 326--339, 1997.

\bibitem{markovsky2021behavioral}
I.~Markovsky and F.~D{\"o}rfler, ``Behavioral systems theory in data-driven
  analysis, signal processing, and control,'' \emph{Ann. Rev. Control},
  vol.~52, pp. 42--64, 2021.

\bibitem{willems2002synthesis}
J.~C. Willems and H.~L. Trentelman, ``{Synthesis of dissipative systems using
  quadratic differential forms: Part I},'' \emph{IEEE Transactions on Automatic
  Control}, vol.~47, no.~1, pp. 53--69, 2002.

\bibitem{willems1986timeI}
J.~C. Willems, ``{From time series to linear system---Part I. Finite
  dimensional linear time invariant systems},'' \emph{Automatica}, vol.~22,
  no.~5, pp. 561--580, 1986.

\bibitem{van2003achievable}
A.~J. Van~der Schaft, ``Achievable behavior of general systems,'' \emph{\SCL},
  vol.~49, no.~2, pp. 141--149, 2003.

\bibitem{julius2005canonical}
A.~A. Julius, J.~C. Willems, M.~N. Belur, and H.~L. Trentelman, ``The canonical
  controllers and regular interconnection,'' \emph{\SCL}, vol.~54, no.~8, pp.
  787--797, 2005.

\bibitem{wolovich1972use}
W.~A. Wolovich, ``The use of state feedback for exact model matching,''
  \emph{\SIAM}, vol.~10, no.~3, pp. 512--523, 1972.

\bibitem{wang1972exact}
S.~Wang and C.~Desoer, ``The exact model matching of linear multivariable
  systems,'' \emph{\TAC}, vol.~17, no.~3, pp. 347--349, 1972.

\bibitem{wang1972solution}
S.~Wang and E.~Davison, ``Solution of the exact model matching problem,''
  \emph{IEEE Transactions on Automatic Control}, vol.~17, no.~4, pp. 574--574,
  1972.

\bibitem{breschi2021direct}
V.~Breschi, C.~De~Persis, S.~Formentin, and P.~Tesi, ``Direct data-driven
  model-reference control with lyapunov stability guarantees,'' in \emph{2021
  60th IEEE Conference on Decision and Control (CDC)}.\hskip 1em plus 0.5em
  minus 0.4em\relax IEEE, 2021, pp. 1456--1461.

\bibitem{willems1991paradigms}
J.~C. Willems, ``Paradigms and puzzles in the theory of dynamical systems,''
  \emph{\TAC}, vol.~36, no.~3, pp. 259--294, 1991.

\bibitem{willems2005note}
J.~C. Willems, P.~Rapisarda, I.~Markovsky, and B.~L.~M. De~Moor, ``A note on
  persistency of excitation,'' \emph{\SCL}, vol.~54, no.~4, pp. 325--329, 2005.

\bibitem{markovsky2020identifiability}
\BIBentryALTinterwordspacing
I.~Markovsky and F.~D{\"o}rfler, ``Identifiability in the behavioral setting,''
  \emph{Vrije Universiteit Brussel, Tech. Rep}, 2020. [Online]. Available:
  \url{{http://homepages.vub.ac.be/~imarkovs/publications/identifiability.pdf.}}
\BIBentrySTDinterwordspacing

\bibitem{francis1976internal}
B.~A. Francis and W.~M. Wonham, ``The internal model principle of control
  theory,'' \emph{Automatica}, vol.~12, no.~5, pp. 457--465, 1976.

\bibitem{markovsky2006exact}
I.~Markovsky, J.~C. Willems, S.~Van~Huffel, and B.~De~Moor, \emph{Exact and
  approximate modeling of linear systems: A behavioral approach}.\hskip 1em
  plus 0.5em minus 0.4em\relax Philadelphia, PA, USA: SIAM, 2006.

\bibitem{coulson2021distributionally}
J.~Coulson, J.~Lygeros, and F.~D{\"o}rfler, ``Distributionally robust chance
  constrained data-enabled predictive control,'' \emph{\TAC}, vol.~67, no.~7,
  pp. 3289--3304, 2021.

\bibitem{van2020willems}
H.~J. van Waarde, C.~De~Persis, M.~K. Camlibel, and P.~Tesi, ``Willems'
  fundamental lemma for state-space systems and its extension to multiple
  datasets,'' \emph{IEEE Control Systems Letters}, vol.~4, no.~3, pp. 602--607,
  2020.

\bibitem{markovsky2022identifiability}
I.~Markovsky and F.~D{\"o}rfler, ``Identifiability in the behavioral setting,''
  \emph{\TAC}, 2022.

\bibitem{padoan2019behavioral}
A.~Padoan, J.~Coulson, H.~J. van Waarde, J.~Lygeros, and F.~D\"orfler,
  ``Behavioral uncertainty quantification for data-driven control,'' in
  \emph{\CDC{61st}}, Cancun, Mexico, 2022, pp. 4726--4731.

\bibitem{yan2023distributed}
Y.~Yan, J.~Bao, and B.~Huang, ``Distributed data-driven predictive control via
  dissipative behavior synthesis,'' \emph{arXiv preprint arXiv:2303.00251},
  2023.

\bibitem{ben2015projectors}
A.~Ben-Israel, ``Projectors on intersection of subspaces,'' \emph{Contemporary
  Mathematics}, vol. 636, pp. 41--50, 2015.

\end{thebibliography}

\end{document}